\documentclass[11pt]{article}
\usepackage{fullpage,epsfig}  
\usepackage{amsmath} 
\usepackage{amssymb}
\usepackage{amsthm}
\usepackage{latexsym}
\usepackage{verbatim}
\usepackage{fixme}
\numberwithin{equation}{section} 
\newtheorem{example}{Example}[section]
\newtheorem{theorem}[example]{Theorem}
\newtheorem{proposition}[example]{Proposition}
\newtheorem{definition}[example]{Definition}
\newtheorem{lemma}[example]{Lemma}

\theoremstyle{definition}
\newtheorem{remark}[example]{Remark}

 \def\norm  #1#2{\Vert\,#1\,\Vert_{#2}}

\newcommand{\R}{{\mathbb R}} 
\newcommand{\CC}{{\mathbb C}} 
\newcommand{\TT}{{\mathcal T}}
\newcommand{\wto}{{\rightharpoonup}}
\newcommand{\N}{{\mathbb N}}
\newcommand{\cA}{\mathcal{A}}
\newcommand{\ups}{\Upsilon^p\mathcal(\R^m)}

\renewcommand{\det}{{\rm det}}

\newcommand{\be}{\begin{eqnarray}}
\newcommand{\ben}{\begin{eqnarray*}}
\newcommand{\een}{\end{eqnarray*}}

\newcommand{\ee}{\end{eqnarray}}

\renewcommand{\d}{{ d}}
\newcommand{\cdm}{{\cal DM}^p_{\cal S}(\O;\R^{m})}

\newcommand{\ra}{\right\rangle}

\newcommand{\la}{\left\langle}
\newcommand{\md}{{\rm d}}

\renewcommand{\O}{\Omega}
\newcommand{\s}{\pi}

\renewcommand{\ker}{{\rm ker}\,{\mathcal A}}
\renewcommand{\b}{\beta}

\newcommand{\Rn}{\R^{n}}

\newcommand{\rca}{\mathcal{M}}

\newcommand{\Hp}{C_{\text{hom}}^p}
\FXRegisterAuthor{jk}{ajk}{JK}

\begin{document}

\begin{sloppypar}

\baselineskip =18pt
\begin{center}
{\huge\bf ${\mathcal A}$-quasiconvexity at the boundary  and weak lower semicontinuity of integral functionals}

\bigskip

{\sc Jan Kr\"{a}mer\footnote{Institute of Mathematics, University of Cologne, 50923 Cologne, Germany} \& Stefan Kr\"{o}mer\footnote{Institute of Mathematics, University of Cologne, 50923 Cologne, Germany}   \& Martin Kru\v{z}\'{\i}k\footnote{Institute of Information Theory and Automation of the CAS,
Pod vod\'{a}renskou v\v{e}\v{z}\'{\i}~4, CZ-182~08~Praha~8, Czech Republic (corresponding address) } \& Gabriel Path\'{o}\footnote{Institute of Information Theory and Automation of the CAS,
Pod vod\'{a}renskou v\v{e}\v{z}\'{\i}~4, CZ-182~08~Praha~8, Czech Republic \& Mathematical Institute, Charles University, Sokolovsk\'{a} 83, CZ-180~00~Praha 8, Czech Republic}
}

\bigskip

\vspace*{.5cm}

\baselineskip=14pt

\vspace*{1.cm}

\begin{minipage}[t]{14cm}

\baselineskip=14pt

{\footnotesize

%\noindent  

{\bf Abstract.}
We state necessary and sufficient conditions for weak lower semicontinuity of $u\mapsto\int_\O h(x,u(x))\,\md x$ where $|h(x,s)|\le C(1+|s|^p)$ is continuous and possesses a recession function, and $u\in L^p(\O;\R^m)$, $p>1$, lives in the kernel of a constant-rank first-order differential operator $\mathcal{A}$ which admits an extension property. Our newly defined notion coincides for $\mathcal{A}=\operatorname{curl}$ with quasiconvexity at the boundary due to J.M.~Ball and J.~Marsden. Moreover, we give an equivalent condition for weak lower semicontinuity of the above functional along sequences weakly converging in $L^p(\O;\R^m)$ and approaching the kernel of $\mathcal{A}$ even if $\mathcal{A}$ does not have the extension property.

\medskip

%\noindent  


\medskip

%\noindent 

}

\end{minipage}

\end{center}
\pagebreak\listoffixmes
\section{Introduction}
In this paper, we 
investigate   the influence of  concentration effects generated by
sequences $\{u_k\}_{k\in\N}\subset L^p(\O;\R^m)$, which   satisfy a linear differential constraint ${\mathcal A}u_k=0$ ($\mathcal{A}$-free sequence), or ${\mathcal A}u_k\to 0$ in $W^{-1,p}(\O;\R^d)$, $1<p<+\infty$ (asymptotically $\mathcal{A}$-free sequence),  on weak lower semicontinuity of integral functionals in the form \fxnote{added since 1409}
\begin{equation}\label{eq:funct}
I(u):=\int_\O h(x,u(x))\,\md x\ .
\end{equation}
 Here, ${\mathcal A}$ is a first-order linear differential operator. To the best of our knowledge, the first such results   were proved  in \cite{fomu} for nonnegative integrands.
In this case, the  crucial necessary and sufficient condition ensuring this property  is the so-called $\mathcal{A}$-quasiconvexity; cf.~\eqref{a-quasiconvexity} below. However, if we refrain from considering only  nonnegative integrands, this condition is not necessarily sufficient.  A prominent example  is 
 $\mathcal{A}$=curl, i.e., when $u$ has a potential.  It is well known  that the weak lower semicontinuity of  $I(u):=\int_\O h(x,u(x))\,\md x$ for $|h(x,s)|\le C(1+|s|^p)$ (i.e. possibly negative and noncoercive)   strongly depends, besides (Morrey's) quasiconvexity, also on the behavior of $h(\cdot,s)$ on the boundary of $\O$. This was first observed by N.~Meyers \cite{meyers} and then elaborated more explicitly  in \cite{kroemer}.   Moreover, it turns out that 
for the special case where $h(x,\cdot)$ possesses a recession function, the precise condition is the  so-called quasiconvexity at the boundary \cite{bama,kruzik}. Namely, if $\{u_k\}_{k\in\N}\subset L^p(\O;\R^m)$ is a weakly converging sequence, concentrations of $\{|u_k|^p\}_{k\in\N}\subset L^1(\O;\R^{m})$ at the boundary of $\O$  can destroy weak lower semicontinuity. We refer to \cite{kkk,kalamajska-kruzik} for a thorough analysis of oscillation and concentration effects in the gradient (curl-free) case.

The situation is considerably more complicated in case of more general operators $\mathcal A$.  In order to see this,  let us isolate a necessary condition for weak lower semicontinuity of $I$ in a simple prototypical situation, a possible candidate to replace quasiconvexity at the boundary for general $\cA$. 
Consider a unit half-ball
$\O:=B(x_0,1)\cap \{x\mid (x-x_0)\cdot\nu_{x_0}\le 0\}\subset\R^n$. We are mainly interested in the behavior near $x_0$, where the boundary of $\O$ is locally flat with normal $\nu_{x_0}$ (a boundary of class $C^1$ actually suffices for the argument below, with some additional technicalities). In addition, we assume for simplicity that the integrand $h=h(u)$ is smooth and positively $p$-homogeneous, i.e.,
for any $\ell\ge 0$, $h(\ell s)=\ell^p h(s)$. Given any $u\in L^p(\R^n;\R^m)\cap\ker$ such that $u$ is compactly supported in $B(0,1)$,
lower semicontinuity along $(u_k)\subset L^p(\R^n;\R^m)\cap\ker$, $u_k(x):=k^{n/p} u(k(x-x_0))$, then implies
$\liminf_{k\to\infty} I(u_k)\ge I(0)=0$, because $u_k\rightharpoonup 0 $ in $L^p$.
Since $I(u_k)=\int_{\O} h(u)\,dx$ for all $k$ by a change of variables, shifting $x_0$ to the origin  
we get a necessary condition on $h$: for all $u\in L^p(B(0,1);\R^m)\cap\ker$ such that $u$ vanishes near the boundary of $B(0,1)$   
\be\label{eq:necessary} 
\begin{aligned}
	&\int_{B(0,1)\cap\{x\cdot\nu_{x_0}\le 0\}}h(u(x))\,\md x\ge 0=  \int_{B(0,1)\cap\{x\cdot\nu_{x_0}\le 0\}}h(0)\,\md x\\
	&\text{for all $u\in L^p(B(0,1);\R^m)\cap\ker$ with $u=0$ near $\partial B(0,1)$}.
\end{aligned}
\ee
%If $\mathcal{A}=$curl, any $u\in L^p(B(0,1);\R^m)\cap\ker$ can be written as the gradient of a potential in $W^{1,p}$, and 
%\eqref{eq:necessary} then becomes
%\be\label{eq:necessarygrad} 
%\begin{aligned}
	%&\int_{B(0,1)\cap\{x\cdot\nu_{x_0}\le 0\}}h(0,\nabla v(x))\,\md x\ge 0
	%&\text{for all $v\in W^{1,p}(B(0,1);\R^k)\cap\ker$ (if $m=kd$) with $v=0$ near $\partial B(0,1)$}.
%\end{aligned}
%\ee
%For the $p$-homogeneous function $h$, this is quasiconvexity at the boundary with respect to the normal $\nu_{x_0}$ (at the zero matrix), in the sense of Ball and Marsden \cite{bama}.
It is clear that for the positively $p$-homogeneous function $h$, this condition generalizes quasiconvexity at the boundary at the zero matrix (for gradients, i.e., curl-free fields) to more general differential constraints given by $\cA$. Together with $\cA$-quasiconvexity, \eqref{eq:necessary} (at every $x_0\in \partial\Omega$, for a smooth domain $\Omega$) is also sufficient for weak lower semicontinuity if $\mathcal{A}=$curl, but as it turns out, this is no longer true for general $\cA$, which also means that \emph{\eqref{eq:necessary} is too weak to act as the correct generalization of quasiconvexity at the boundary} for our purposes:
\begin{example}
Let $n=m=2$, $p=2$. We take $\cA$ to be the differential operator of the Cauchy-Riemann system, i.e., $\mathcal{A}u=0$ if and only if $\partial_1 u_1-\partial_2 u_2=0=\partial_2u_1+\partial_1u_2$ (which in turn means that $u_1+iu_2$ is holomorphic on its domain, as a function of $z=x_1+ix_2\in \CC$). Then \eqref{eq:necessary} is trivially satisfied for any function $h$ with $h(0)=0$, because the only admissible $u$ is the zero function. 
Similarly, any $h$ is $\cA$-quasiconvex, as $\cA$-quasiconvexity is tested with periodic functions in $\ker$ with zero mean, and for Cauchy-Riemann, the only such function is zero due to the Liouville theorem. 
Nevertheless, for $h(x,s):=-|s|^2$ and any bounded domain $\O\subset\R^2\cong \CC$ with smooth boundary, $I$ is \emph{not} weakly lower semicontinuous in $L^p\cap \ker$: Let $u_k(z)=\frac{1}{k(z-z_k)}$, where $\{z_k\}\subset \CC \setminus\O$ is a sequence defined in such a way that 
$\int_\O |u_k(x)|^2\,\md x=1$ (there always exists such a $z_k$ by continuity, because for fixed $k$, $\int_\O |u_k|^2\,dx\to 0$ as $|z_k|\to \infty$ and $\int_\O |u_k|^2\,dx\to +\infty$ as 
$\operatorname{dist}(z_k;\Omega)\to 0$). In particular, $z_k$ approaches the boundary of $\O$ from the outside as $k$ increases. Then $u_k\rightharpoonup 0$ in $L^2(\O;\R^2)$ but $\liminf_{k\to\infty} I(u_k)=-1<I(0)=0$. 
\end{example}
The example shows that test functions in the operator kernel and with zero `` boundary conditions''  do not suffice  to  analyze concentration effects on the boundary like that of our holomorphic sequence $u_k$ in the example, where a singularity is approaching the boundary from the outside. Replacing the class of test functions in \eqref{eq:necessary} by periodic functions with zero mean as in the definition of $\cA$-quasiconvexity does not help either, because \eqref{eq:necessary}  would still be trivially satisfied in the example, now due to the Liouville theorem. Altogether, we see that the problem of weak lower semicontinuity for a generic $\mathcal{A}$ is considerably more involved, once negative integrands are allowed.

Nevertheless, sequences of functions with smaller and smaller support are certainly natural to test weak lower semicontinuity along ``point concentrations''. The only question is how that should be reflected in an appropriate stronger version of \eqref{eq:necessary}.
This dilemma is resolved below in Definitions~\ref{def:aaqcb} and \ref{def:aqcb} by 
allowing test functions to depart (in a controlled way)  from the kernel of $\mathcal{A}$. We show that this approach naturally gives a new necessary and sufficient condition for weak lower semicontinuity of $I$ along asymptotically $\mathcal{A}$-free sequences (${\mathcal A}u_k\to 0$) called here {\it strong $\mathcal{A}$-quasiconvexity at the boundary}; cf.~Def.~\ref{def:aaqcb}, even for quite rough domains. Obviously, strong   $\mathcal{A}$-quasiconvexity at the boundary also suffices for wlsc of $I$ along sequences in the kernel of $\mathcal{A}$. We also derive a necessary and sufficient condition for the latter situation, called {\it $\mathcal{A}$-quasiconvexity at the boundary};  cf.~Def.~\ref{def:aqcb}. As the name suggests, strong  $\mathcal{A}$-quasiconvexity at the boundary implies $\mathcal{A}$-quasiconvexity at the boundary, but in general, these notions are not equivalent as outlined in Section~\ref{sec:concluding}, where we also discuss a sufficient condition on the operator $\mathcal{A}$ and the domain ensuring equivalence (Def.~\ref{def:afe}).
The picture is therefore more complicated than in the case of nonnegative integrands $h$, where weak $\mathcal{A}$-quasiconvexity (see Def.~\ref{a-quasiconvexity})  of $h(x,\cdot)$ is known to be a necessary and sufficient condition for weak lower semicontinuity of $I$ \cite[Thm~3.6, 3.7]{fomu} in both cases, i.e. if  ${\mathcal A}u_k=0$ or if ${\mathcal A}u_k\to 0$  in $W^{-1,p}(\O;\R^d)$.

%The notion of $\mathcal{A}$-quasiconvexity at the boundary and  strong $\mathcal{A}$-quasiconvexity at the boundary are first defined for general domains 
%in Definitions~\ref{def:aqcb} and \ref{def:aaqcb}, respectively and then we derive equivalent simpler  versions for domains with $C^1$ boundaries. See Prop.~\ref{prop:aAqcbD}, Prop.~\ref{prop:xAqcbD}, and Prop.~\ref{prop:aqcb}. However, even for Lipschitz domains such a simplification is not really possible.  

Let us emphasize  that variational problems with differential constraints naturally appear in hyperelasticity, electromagnetism, or in micromagnetics \cite{desimone,pedregal0,pedregal} and are closely related to the theory of compensated compactness \cite{murat1,tartar,tartar1}. 
The concept of $\mathcal{A}$-quasiconvexity  goes back to \cite{dacorogna82} and  has been proved to be useful as a unified approach to variational problems with differential constraints, including results on homogenization \cite{braides,fonseca-kroemer}, dimension reduction \cite{kreisbeck-rindler} and characterization of  generalized Young measures \cite{baia-matias-santos} in the $\mathcal{A}$-free setting. Moreover, first results on $\mathcal{A}$-quasiaffine functions and weak continuity appeared in \cite{foss-rand}.
As to weak lower semicontinuity,  the theory was first developed for nonnegative integrands in \cite{fomu} as mentioned before, with  extensions to nonnegative functionals with nonstandard growth \cite{fonseca-leoni-mueller} and the case of an operator $\mathcal{A}$ with nonconstant coefficients \cite{santos}. The recent work \cite{baia} analyzes lower semicontinuity 
of functionals with linearly growing integrands, including negative integrands but excluding concentrations at the domain boundary.

The plan of the paper is as follows.  We first recall some needed definitions and results in Section~\ref{preliminaries}. Our newly derived conditions which, together 
with $\mathcal{A}$-quasiconvexity  precisely characterize weak lower semicontinuity are studied in Section~\ref{sec:aqcb}. The main results are summarized in Theorem~\ref{thm:wlscasympt} and Theorem~\ref{thm:wlscAfree}. After the concluding remarks in the final section, some auxiliary material is provided in the appendix.

\bigskip

\section{Preliminaries}\label{preliminaries}

Unless explicitly stated otherwise, we always work with a bounded  domain  $\O\subset\R^n$ such that $\mathcal{L}^n(\partial\O)=0$, equipped with the Euclidean topology and the $n$-dimensional Lebesgue measure $\mathcal{L}^n$. 
 $L^p(\O,\R^m)$, $1\le p\le +\infty$,  is a standard  Lebesgue space.   Furthermore, $W^{1,p}(\O;\R^m)$, $1\le p<+\infty$,  stands for  the usual space of measurable mappings, which together with 
their first (distributional) derivatives, are integrable with the $p$-th power. A space of mappings from   $W^{1,p}(\O;\R^m)$ with zero traces  is  standardly    denoted  $W^{1,p}_0(\O;\R^m)$.
If $1<p<+\infty$   then $W^{-1,p}(\O;\R^m)$  denotes the dual space to $W_0^{1,p'}(\O;\R^m)$, where  $p'^{-1}+p^{-1}=1$. A sequence $\{u_k\}_{k\in\N}$ converges to zero in measure if $\mathcal{L}^n(\{x\in\O;\,|u_k(x)|\ge\delta\})\to 0$ as $k\to\infty$, for every $\delta>0$. 

We say that $v\in\ups$ if  there exists a continuous and positively $p$-homogeneous function $v_\infty:\R^m\to\R$, i.e.,  $v_\infty(t s)=t^p v_\infty(s)$ for all $t\ge 0 $ and $ s\in\R^m$, such that 
\be\label{recessionf}
\lim_{|s|\to\infty}\frac{v(s)-v_\infty(s)}{|s|^p}=0\ .
\ee 
 
Such a function is called the {\it recession function} of $v$.

\subsection{The operator ${\mathcal A}$ and  ${\mathcal A}$-quasiconvexity}
Following \cite{fomu},  we consider linear operators $A^{(i)}:\R^m\to\R^d$, $i=1,\ldots, n$, and  define ${\mathcal A}:L^p(\O;\R^m)\to W^{-1,p}(\O;\R^d)$ by
\ben
{\mathcal A}u:=\sum_{i=1}^nA^{(i)}\frac{\partial u}{\partial x_i}\ ,\mbox{where }\ u:\O\to\R^m\ ,
\een
i.e.,
for all $w\in W^{1,p'}_0(\O;\R^d)$ 
\begin{eqnarray*}
\la {\mathcal A}u,w\ra=-\sum_{i=1}^n\int_\O A^{(i)}u(x)\cdot\frac{\partial w(x)}{\partial x_i}\,\md x
\ .\end{eqnarray*}

For   $w\in\R^n$ we define the linear map 
\ben
{\mathbb A}(w):=\sum_{i=1}^n w_iA^{(i)}\ :\ \R^m\to\R^d\ .\een
Throughout this article, we assume that there is $r\in\N\cup\{0\}$ such that
\be\label{rank}
{\rm rank}\ {\mathbb A}(w)=r\mbox{ for all $w\in\R^n$}\ , |w|=1\ ,\ee
i.e.,  ${\mathcal A}$ has the so-called {\it constant-rank property}.

Below, we use $\ker$ to denote the set of all locally integrable functions $u$ such that
$\cA u=0$ in the sense of distributions, i.e., $\int u\cdot \cA^* w\,dx=0$ for all $w\in C^\infty$ compactly supported in the domain, where $\cA^*$ is the formal adjoint of $\cA$. 
Of course, this depends on the domain considered, which should be clear from the context.
In particular, a periodic function $u$
in the space
\[
	L^p_\#(\R^n;\R^m):=\{u\in L^p_{\rm loc}(\R^n;\R^m): \mbox{ $u$ is $Q$-periodic}\}\ 
\]
is in $\ker$ if and only if $\cA u=0$ on $\R^n$. Here and in the following,
$Q$ denotes the unit cube $(-1/2,1/2)^n$ in $\R^n$, and we say that $u:\R^n\to\R^m$ is {\it $Q$-periodic} if for all $x\in\R^n$ and all 
$z\in\mathbb{Z}^n$ 
$$
u(x+z)=u(x)\ .$$

We will use the following lemmas proved in \cite[Lemma 2.14]{fomu} and \cite[Lemma~2.15]{fomu}, respectively.

\begin{lemma}(projection onto $\cA$-free fields in the periodic setting)\label{T}
There is a linear bounded operator $\TT:L^p_\#(\R^n;\R^m)\to L^p_\#(\R^n;\R^m)$ that vanishes on constant functions, $\TT(\TT u)=\TT u$ for all $u\in L^p_\#(\R^n;\R^m)$, and $\TT u\in\ker$. Moreover, for all  $ u\in L^p_\#(\R^n;\R^m)$ with $\int_{Q} u(x)\,\md x=0$ it holds that  
$$
\|u-\TT u\|_{L^p_\#(\R^n;\R^m)}\le C\|{\mathcal A}u\|_{W^{-1,p}_\#(\R^n;\R^d)}\ ,
$$
where $C>0$ is a constant independent of $u$ and $W^{-1,p}_\#$ denotes the dual of $W^{1,p'}_\#$ ($\frac{1}{p'}+\frac{1}{p}=1$), 
the $Q$-periodic functions in $W^{1,p'}_{\text{loc}}(\R^n;\R^m)$ equipped with the norm of $W^{1,p'}(Q;\R^m)$.
\end{lemma}

\begin{remark}\label{rem:T}
For every $w\in W^{-1,p}_\#(\R^n)$, we have $\|w\|_{W^{-1,p}(Q)}\leq \|w\|_{W^{-1,p}_\#(\R^n)}$. The converse inequality does not hold, not even up to a constant. However, Lemma~\ref{T} is often applied to (a sequence of) functions supported in a fixed set $G\subset\subset Q$ (up to periodicity, of course). One can always find a constant $C=C(\O,p,G)$ such that
\[
	\|\cA u\|_{W^{-1,p}_\#(\R^n;\R^d)} \leq C \|\cA u\|_{W^{-1,p}(Q;\R^d)}\quad\text{for every $u\in L^p(Q;\R^m)$ with $u=0$ a.e.~on $Q\setminus G$.}
\]
To achieve this, the $Q$-periodic test functions used in the definition of the norm in $W^{-1,p}_\#$ can be multiplied with a fixed cut-off function $\eta\in C_0^\infty(Q;[0,1])$ with $\eta=1$ on $G$ to make them admissible (i.e., elements of $W_0^{1,p'}(Q)$) for the supremum defining the norm in $W^{-1,p}$. This enlarges their norm in $W^{1,p'}$ at most by a constant factor which only depends on $p$ and $\|\nabla \eta\|_{L^\infty(Q)}$ (and thus the distance of $G$ to $\partial Q$).
\end{remark}

\begin{lemma}(Decomposition Lemma)\label{decomposition}
Let $\O\subset\R^n$ be bounded and open, $1<p<+\infty$, and let $\{u_k\}\subset L^p(\O;\R^m)$ be bounded and such that  ${\mathcal A}u_k\to 0$ in $W^{-1,p}(\O;\R^d)$ strongly, and $u_k\wto u$ in $L^p(\O;\R^m)$ weakly. Then there is a sequence $\{z_k\}_{k\in\N}\subset L^p(\Omega;\R^m)\cap\ker $, $\{|z_k|^p\}$ is equiintegrable in $L^1(\Omega)$ and $u_k-z_k\to 0$ in measure in $\Omega$.
\end{lemma}

We also point out the following simple observation made in the proof of Lemma 2.15 in \cite{fomu}, which is useful if we need to truncate
$\cA$-free or ``asymptotically'' $\cA$-free sequences:
\begin{lemma}\label{lem:truncate}
Let $\O\subset\R^n$ be open and bounded, and let $\{u_k\}\subset L^p(\O;\R^m)$ be a bounded sequence such that ${\mathcal A}u_k\to 0$ in $W^{-1,p}(\O;\R^d)$ strongly and $u_k\wto 0$ in $L^p(\O;\R^m)$ weakly. Then for every $\eta\in C^\infty(\R^n)$, ${\mathcal A}(\eta u_k)\to 0$ in $W^{-1,p}(\O;\R^d)$.
\end{lemma}
{\it Proof.}
$\mathcal{A}(\eta u_k)=\eta\mathcal{A}u_k+\sum_{i=1}^n u_kA^{(i)}\frac{\partial\eta}{\partial x_i}\to 0$ in $W^{-1,p}$, the second term due to the compact embedding of $L^p$ into $W^{-1,p}$.
\hfill $\Box$

\begin{definition}\label{a-quasiconvexity} (see \cite[Def.~3.1,~3.2]{fomu})
 We say that  a continuous function $v:\R^m\to\R$, $|v|\le C(1+|\cdot|^p)$ for some $C>0$, is ${\mathcal A}$-{\it quasiconvex} if for all $s_0\in\R^{m}$ and all $\varphi\in L^p_\#(Q; \R^m)\cap\ker$ with $\int_{Q}\varphi(x)\,\md x=0$ it holds
\ben
v(s_0)\le \int_Q v(s_0+\varphi(x))\,\md x\ .\een
\end{definition}

Besides curl-free fields, admissible examples of $\mathcal {A}$-free mappings  include solenoidal fields where $\mathcal{A}=\rm{div}$ and higher-order gradients where $\mathcal {A}u=0$ if and only if $u=\nabla^{(s)}\varphi$ for some $\varphi\in W^{s,p}(\O;\R^\ell)$, and some $s\in\N$ (for more details see Subsection~\ref{subsec:hessian}, where $s=2$).

\subsection{Weak lower semicontinuity}

Let $I:L^p(\O;\R^m)\to\R$ be defined as
\be \label{def:funI}
I(u):=\int_\O h(x,u(x))\,\md x\ .
\ee
%We often restrict $I$ to $\ker$ below.

\begin{definition}\label{def:asymptAfree} \hfill

\noindent (i) We say that a sequence $\{u_k\}\in L^p(\O;\R^m)$ is asymptotically $\cA$-free if $\|\cA u_k\|_{W^{-1,p}(\O;\R^m)} \to 0$ as $k\to \infty$. 

\noindent (ii) A functional $I$ as in \eqref{def:funI} is called weakly sequentially lower semicontinuous (wslsc) along asymptotically $\cA$-free sequences in $L^p(\O;\R^m)$ if $\liminf_{k\to\infty} I(u_k)\geq I(u)$ for all such sequences that weakly converge to some limit $u$ in $L^p$.

\noindent (iii) Analogously, we say that   a functional $I$ is weakly sequentially lower semicontinuous (wslsc) along $\mathcal{A}$-free sequences in  $L^p(\O;\R^m)$ if \[\liminf_{k\to\infty} I(u_k)\geq I(u)\text{ for  all }\{u_k\}\subset L^p(\O;\R^m)\cap\ker.\] 
\end{definition}

\bigskip

We have the following result which was proved in \cite[Theorem 2.4]{ifmk} in a slightly less general version. However, its original proof directly extends to this setting. 

\begin{theorem}\label{wlsc1}
Let $h:\bar\O\times\R^m\to\R$ be continuous such that $h(x,\cdot)\in\ups$ for all $x\in\bar\O$ and $h(x,\cdot)$ is $\mathcal{A}$-quasiconvex for almost every $x\in\O$, $1<p<+\infty$. Then  $I$  is sequentially weakly
lower semicontinuous in $L^{p}(\O;\R^m)\cap\ker $ if and only if for any bounded sequence $\{u_k\}\subset L^{p}(\O;\R^m)\cap\ker $ such that $u_k\to 0$ in measure there is
\be\label{wlsc-inequality}
\liminf_{k\to\infty}I(u_k)\ge I(0)\ .
\ee
\end{theorem}
The statement of Theorem~\ref{wlsc1} remains valid if we replace the sequences in $\ker$ with asymptotically $\cA$-free sequences. 

\begin{theorem}\label{thm:wlscaAfree}
With $h$ and $p$ as in Theorem~\ref{wlsc1}, $I$ is wslsc along asymptotically $\cA$-free sequences in $L^{p}(\O;\R^m)$ if and only if \eqref{wlsc-inequality} holds
for any
bounded, asymptotically $\cA$-free sequence $\{u_k\}\subset L^{p}(\O;\R^m)$ such that $u_k\to 0$ in measure.
\end{theorem}
{\it Proof.}
We only point out the differences to the proof \cite[Theorem 2.4]{ifmk}. First, the result there is stated only for functions $h$ of product form $h(x,\xi)=g(x)v(\xi)$, but as in the case of Theorem~\ref{wlsc1}, it works verbatim also for our slightly more general class. ``Only if'' is trivial as before. For ``if'', we also rely on splitting a given sequence into a purely oscillating ($p$-equiintegrable) part and a purely concentrating part, which is still a straightforward application of the decomposition lemma (Lemma~\ref{decomposition}). Notice that the purely oscillating part $\{z_k\}$ lives in $\ker$, even if the sequence we started with is only asymptotically $\cA$-free. The rest of the proof is completely analogous to the corresponding one in \cite{ifmk}.
\hfill$\Box$ 

\bigskip

\begin{remark}\label{replacement} \hfill

\noindent (i) It follows from \cite[(5.1)]{ifmk} that \eqref{wlsc-inequality} can be replaced by 
\[ 
	\liminf_{k\to\infty}I_\infty(u_k)\ge I_\infty(0)=0,~~\text{where}~I_\infty(u):=\int_\O h_\infty(x,(u(x))\,\md x,
\] 
with $h_\infty(x,\cdot)$ denoting the recession function of $h(x,\cdot)$.

\noindent (ii) In fact, having an integrand $(x,s)\mapsto h(x,s)$ which is $\mathcal{A}$-quasiconvex in the second variable, weak lower semicontinuity can only fail due to sequences concentrating large values on small sets, and it even suffices to test that with sequences $\{u_k\}$ which tend to zero in measure and concentrate at the boundary in the sense that $\{|u_k|^p\}$ converges weakly* to a measure $\sigma\in\mathcal{M}(\bar\O)$ with $\sigma(\partial\O)>0$.   
\end{remark}

\bigskip

\section{Notions of $\mathcal{A}$-quasiconvexity at the boundary}\label{sec:aqcb}

The two conditions introduced below play a crucial role in our characterization of weak lower semicontinuity of integral functionals. They are typically applied to the recession function $h_\infty$ of an integrand $h$ with $p$-growth.

Before we state them, we fix some additional notation frequently used in what follows:
\[
\begin{aligned}
	&L^p_0(\O;\R^m):=\{u\in L^p(\O;\R^m);\, \mbox{supp}\, u\subset\O\},\\
	&\Hp(\R^m):=\{v\in C(\R^m);\, v\mbox{ is positively $p$-homogeneous}\}.
\end{aligned}
\]
A norm in $\Hp$ is given by the supremum norm taken on the unit sphere in $\R^m$. Moreover, whenever a larger domain comes into play, functions in $L^p_0(\O;\R^m)$ are understood to be extended by zero to $\R^n\setminus \O$ without changing notation.

The definitions given below are stated in a form which is suitable for rather general domains and the most natural in the proofs of our characterizations of weak lower semicontinuity. For domains with a boundary of class $C^1$, equivalent, simpler variants more closely resembling the original notion of quasiconvexity at the boundary in the sense of Ball and Marsden are presented in Proposition~\ref{prop:aAqcbD}--Proposition~\ref{prop:aqcb}.
\begin{definition}\label{def:aqcb}
We say that $h_\infty\in C(\bar\O;\Hp(\R^m))$ is $\cA$-quasiconvex at the boundary ($\mathcal{A}$-qcb) at $x_0\in\partial\O$ if for every  $\varepsilon>0$ there are $ \delta>0$ and $\alpha>0$  such that 
\be\label{a-qcb} 
	\int_{B(x_0,\delta)\cap\Omega}  h_\infty(x,u(x))+\varepsilon|u(x)|^p\,\md x\ge 0\ 
\ee
for every $u\in L^p_0(B(x_0,\delta);\R^m)$ with $\|\mathcal{A}u\|_{W^{-1,p}(\R^n;\R^d)}<\alpha \|u\|_{L^p(B(x_0,\delta)\cap\O;\R^m)}$.
\end{definition}

The next notion is intimately related to weak lower semicontinuity along asymptotically $\mathcal{A}$-free sequences. Notice that the {\it  only  but crucial}  difference between Definitions~\ref{def:aqcb} and \ref{def:aaqcb}  is the norm used to measure $\mathcal{A}u$. 

\begin{definition}\label{def:aaqcb}
We say that $h_\infty\in C(\bar\O;\Hp(\R^m))$ is  strongly $\cA$-quasiconvex at the boundary (strongly-$\mathcal{A}$-qcb) at $x_0\in\partial\O$ if for every  $\varepsilon>0$ there are $ \delta>0$ and $\alpha>0$  such that 
\be\label{aa-qcb} 
	\int_{B(x_0,\delta)\cap\Omega}  h_\infty(x,u(x))+\varepsilon|u(x)|^p\,\md x\ge 0\ 
\ee
for every $u\in L^p_0(B(x_0,\delta);\R^m)$ with $\|\mathcal{A}u\|_{W^{-1,p}(\O;\R^d)}<\alpha \|u\|_{L^p(B(x_0,\delta)\cap\O ;\R^m)}$.\jknote{Intersection with $B(x_0,\delta)$ made explicit.}
\end{definition}

As it turns out, strong $\mathcal{A}$-qcb is natural in the characterization for weak lower semicontinuity along asymptotically $\cA$-free sequences, while $\mathcal{A}$-qcb plays the same role for weak lower semicontinuity along precisely $\cA$-free sequences. While strong $\mathcal{A}$-qcb always implies $\mathcal{A}$-qcb, they are not equivalent in general (see Section~\ref{sec:concluding}).

\begin{remark}\label{rem:Aqcb-extension}
Due to the fact that $\mathcal{A}u$ in Definition~\ref{def:aqcb} is required to be small on $B(x_0,\delta)$, a set which is not fully contained in $\O$, $\mathcal{A}$-qcb as defined above can only natural if there is an $\cA$-free extension operator on $L^p(\Omega;\R^m)$, cf.~Definition~\ref{def:Aext} below. However, the existence of such an extension operator may require sufficient smoothness of $\partial \O$, and, worse, it strongly depends on $\cA$ (it fails for the Cauchy-Riemann system, e.g.). 
The strong variant of $\mathcal{A}$-qcb does not have this unpleasant implicit dependence on $\cA$-free extension properties.
\end{remark}

\begin{remark}\label{rem:Aqcb-test}
In Definition~\ref{def:aqcb}, $\cA u$ is measured in the norm of $W^{-1,p}(\R^n;\R^d)$, but $\R^n$ can be replaced by any domain $S_\delta$ compactly containing $B(x_0,\delta)$, because for distributions supported on $B(x_0,\delta)$, the norms of $W^{-1,p}(\R^n;\R^d)$ and $W^{-1,p}(S_\delta;\R^d)$ are equivalent, with constants depending on $\delta$. The latter is not a problem since $\alpha$ depends on $\varepsilon$ and thus may also depend on $\delta=\delta(\varepsilon)$. In particular, $\mathcal{A}$-qcb can also be defined using the class of all $u\in L^p_0(B(x_0,\frac{\delta}{2});\R^m)$ with $\|\mathcal{A}u\|_{W^{-1,p}(B(x_0,\delta);\R^d)}<\alpha \|u\|_{L^p(B(x_0,\delta)\cap\O;\R^m)}$.
Similarly, the class of test functions in Definition~\ref{def:aaqcb}  can be replaced by the set of all $u\in L^p_0(B(x_0,\frac{\delta}{2});\R^m)$ such that $\|\mathcal{A}u\|_{W^{-1,p}(\O\cap B(x_0,\delta);\R^d)}<\alpha \|u\|_{L^p(\O\cap B(x_0,\delta);\R^m)}$.
\end{remark}

%Notice that the smallness of $\cA u$ is  measured in the $W^{-1,p}$-norm on $\O$ instead of a larger set as in Definition~\ref{def:aqcb}. To calculate this norm, we seek the largest possible value of $\int_{\R^n} u\cdot \cA^* \varphi \md x$ among all functions  $\varphi \in W_0^{1,p'}(\O;\R^d)$ with norm not larger than $1$ in that space. In particular, each admissible $\varphi$ is now required to vanish on $\partial \O$. This {\it does} make a difference, which can be easily be checked in a simplified setting: if we let $B$ denote a ball in $\R^n$ and $D=\{x\in B|x\cdot \nu<0\}$ the half ball in $B$ determined by some (arbitrary but fixed) vector $\nu$, then the norms of the dual spaces of $X:=W_0^{1,p'}(D)$ and $Y:=\{v\in W_0^{1,p'}(D)\mid v=u|_D~\text{for some}~u\in W_0^{1,p'}(B)\}$ are \emph{not equivalent} on their mutual subset $L^p_0(B)$.

\begin{remark}
In Definition~\ref{def:aqcb} as well as in Definition~\ref{def:aaqcb}, if for a given $\varepsilon>0$ the estimate holds for some $\delta>0$, then it also holds for any $\tilde\delta<\delta$ in place of $\delta$. Hence, both $\cA$-qcb and strong $\cA$-qcb are local properties of $h_\infty$ in the $x$ variable, since it suffices to study arbitrarily small neighborhoods of $x_0$.
\end{remark}

It is possible to formulate several equivalent variants of the definitions of $\cA$-quasiconvexity at the boundary. In particular, the following proposition shows that the first variable of $h$ can be ``frozen'' in Definition~\ref{def:aaqcb}.  

\begin{proposition}~\label{prop:sAqcb-frozen}
 A function $(x,s)\mapsto h_\infty(x,s)$, $h_\infty\in C(\bar\O;\Hp(\R^m))$, is strongly $\mathcal{A}$-qcb at $x_0\in\partial\Omega$ if and only if $s\mapsto h_\infty(x_0,s)$ is strongly $\mathcal{A}$-qcb at $x_0\in\partial\Omega$.
\end{proposition}

\begin{comment}
A function $h_\infty\in C(\bar\O;\Hp(\R^m))$,
%$h:\bar\O\times\R^m\to\R$ such that $h_\infty(\cdot,s)\in C(\bar\O)$ and $h(x,\cdot)$ positively $p$-homogeneous for all $x\in\bar\O$
is $\mathcal{A}$-qcb at $x_0\in\partial\O$ if and only if 
 for all  $\varepsilon>0$ there are $ \delta>0$, $\alpha>0$  such that 
 for all $u\in L^p_0(B(x_0,\delta);\R^m)$ with $\|\mathcal{A}u\|_{W^{-1,p}}<\alpha \|u\|_{L^p}$ (all norms taken on $B(x_0,\delta)\cap \O$),
 \be\label{a-qcb0} \int_{B(x_0,\delta)\cap\Omega}  h_\infty(x_0,u(x))+\varepsilon|u(x)|^p\,\md x\ge 0\ .\ee
\end{comment}

{\it Proof.}
Let $\varepsilon>0$ and recall that if \eqref{a-qcb} holds for some $\delta>0$ then it holds also for any $0<\tilde\delta<\delta$ in the place of $\delta$.
We have 
\[
\begin{aligned}
\left|\int_{B(x_0,\delta)\cap\O} h_\infty\left(x,\frac{u(x)}{|u(x)|}\right)|u(x)|^p\,\md x- \int_{B(x_0,\delta)\cap\O} h_\infty\left(x_0,\frac{u(x)}{|u(x)|}\right)|u(x)|^p\,\md x\right|&\\
\le \int_{B(x_0,\delta)\cap\O}\mu(|x-x_0|,0)|u(x)|^p\,\md x\le M(\delta)\int_{B(x_0,\delta)\cap\O}|u(x)|^p\,\md x\ &,
\end{aligned}
\]
where $\mu:\R\times\R\to\R$ is a continuous modulus of continuity of  the continuous function $h_\infty$ restricted to the compact set $\bar\O\times S^{m-1}$ and $M(\delta):=\max_{x\in \overline{B(x_0,\delta)\cap\O}}\mu(|x-x_0|,0)$. In particular, $M(\delta)\to 0$ as $\delta\to 0$. Hence, if \eqref{a-qcb} holds then we have that 
$$
\int_{B(x_0,\delta)\cap\Omega} h_\infty(x_0,u(x))+(M(\delta)+\varepsilon)|u(x)|^p\,\md x\ge \int_{B(x_0,\delta)\cap\Omega} h_\infty(x,u(x))+\varepsilon|u(x)|^p\,\md x\ge 0\ ,
$$
whence $(x,s)\mapsto h_{\infty}(x_0,s)$ is strongly $\cA$-qcb at $x_0$. Here, note that $M(\delta)+\varepsilon$ can be made arbitrarily small if $\delta$ is small enough. The converse implication is proved analogously.
\hfill $\Box$

\bigskip

Exactly as in the case of Definition~\ref{def:aaqcb}, the first variable of $h_\infty$ can be ``frozen'' in Definition~\ref{def:aqcb}:

\begin{proposition}
 A function $(x,s)\mapsto h_\infty(x,s)$, $h_\infty\in C(\bar\O;\Hp(\R^m))$, is $\mathcal{A}$-qcb at $x_0\in\partial\Omega$ if and only if 
$(x,s)\mapsto h(x_0,s)$ is $\mathcal{A}$-qcb at $x_0\in\partial\Omega$.
\end{proposition}
\begin{comment}
A function $h_\infty\in C(\bar\O;\Hp(\R^m))$ is $\mathcal{A}$-qcb at $x_0\in\partial\O$ if and only if for all  $\varepsilon>0$ there are $ \delta>0$, $\alpha>0$  such that for all $u\in L^p_0(B(x_0,\delta);\R^m)$ with $\|\mathcal{A}u\|_{W^{-1,p}(\R^n;\R^d)}<\alpha \|u\|_{L^p(\O;\R^m)}$,
 \be\label{a-qcb0B} 
\int_{B(x_0,\delta)\cap\Omega}  h_\infty(x_0,u(x))+\varepsilon|u(x)|^p\,\md x\ge 0\ .
\ee
\end{comment}

By itself, ``freezing'' the first variable of $h$ does not help to really simplify Definition~\ref{def:aqcb} or Definition~\ref{def:aaqcb},
because the possibly complicated local shape of the domain can still prevent passing to the limit as $\delta\to 0$ in a meaningful way.
However, this is the best we can do without imposing further restrictions on the smoothness of $\partial\Omega$. Even for Lipschitz domains, the general form of the Definitions typically cannot be avoided (see \cite[Remark 1.8]{kroemer} for a more detailed discussion of this in the gradient case corresponding to $\cA$=curl).

So far, is not clear to what extent the notion of (strong) $\mathcal{A}$-qcb depends on the local shape of $\partial\Omega$ near the boundary point under consideration. The propositions below show that at least for domains with smooth boundary we can in some sense pass to the limit as $\delta\to 0$
in Definition~\ref{def:aqcb} and Definition~\ref{def:aaqcb}, and the domain enters only via the outer normal to $\partial\O$ at this point.

\begin{proposition}\label{prop:aAqcbD}
Assume that $\O\subset\R^n$ has a $C^1$-boundary in a neighborhood of $x_0\in\partial\O$. Let $\nu_{x_0}$ be the outer unit normal to $\partial\O$ at $x_0$ and
\[
	D_{x_0}:=\{x\in B(0,1)\mid \,x\cdot\nu_{x_0}< 0\}.
\]
Then  $v\in \Hp(\R^m)$ is  strongly-$\mathcal{A}$-qcb at $x_0$ if and only if
\be\label{half-ball}
\begin{aligned}
	&\text{for every $\varepsilon>0$ there exists $\beta>0$ such that}\\
 	&\qquad\int_{D_{x_0}} v(\varphi(x))+\varepsilon|\varphi(x)|^p\,\md x\ge 0\   \\
	&\text{for every $\varphi\in L^p_0(B(0,\tfrac{1}{2});\R^m)$ with
	$\|\mathcal{A}\varphi\|_{W^{-1,p}(D_{x_0};\R^d)}\le\beta \|\varphi\|_{L^p(D_{x_0};\R^m)}$.}
\end{aligned}
\ee
\end{proposition}

{\it Proof.} Without loss of generality let us assume $x_0=0$.  We adapt the proof which appeared already in \cite{kroemer} for the gradient case.

``only if'': Suppose that $v$ is strongly-$\mathcal{A}$-qcb at $0$.  Take $\varepsilon>0$ and get $\alpha$, $\delta>0$ such that
\be \label{ineq_u}
\int_{B(0,\delta) \cap \O}v(u(x))+\varepsilon|u(x)|^p\,\md x\ge 0
\ee
for every $u\in L^p_0(B(0,\tfrac{\delta}{2});\R^m)$ satisfying $\|\mathcal{A}u\|_{W^{-1,p}(B(0,\delta) \cap \O;\R^d)} \le \alpha \|u\|_{L^p(B(0,\delta) \cap \O;\R^m)}$.  Introducing the scaling $\Phi_\delta \colon \overline{B(0,\delta)} \ni x \mapsto \delta^{-1} x \in \overline{B(0,1)}$, the inequality \eqref{ineq_u} can be rewritten as
\be \label{ineq_y}
\int_{\delta^{-1}(\O\cap B(0,\delta))}v(y(x'))+\varepsilon|y(x')|^p\,\md x'\ge 0\ , \text{ where $y := \delta^{n/p} u \circ \Phi_\delta^{-1}$} \label{ineq64261}
\ee
Due to the smoothness of the boundary near zero, there exists a transformation $\Psi_\delta \colon \overline{B(0,1)}\to\overline{B(0,1)}$ such that $\Psi_\delta(0)=0$, $\Psi_\delta(B(0,\tfrac{1}{2}))=B(0,\tfrac{1}{2})$ and $\Psi_\delta(D_{0})=\delta^{-1}(\O \cap B(0,\delta))$, while both $\Psi_\delta$ and its inverse $\Psi^{-1}_\delta$ converge to the identity in $C^1(\overline{B(0,1)};\R^n)$ as $\delta\to 0$. Hence, \eqref{ineq_y} leads to
\be \label{ineq_phi}
\int_{D_{0}}(v(\varphi(z))+\varepsilon|\varphi(z)|^p)|\det {\rm D}_z \Psi_\delta(z)|\,\md z\ge 0\ ,
\ee
where $\varphi := y \circ \Psi_\delta$ and $[{\rm D}_z \Psi_\delta]_{ij} := \partial \Psi_{\delta i} / \partial z_j$ for $i$, $j = 1, \dots, n$.  Due to the boundedness of $v + \varepsilon|\cdot|^p$  and the (uniform) continuity of the transformation $\Psi_\delta$ on the unit sphere, we have the estimate
\be \label{ineq_estim}
|(v(\varphi(z))+\varepsilon|\varphi(z)|^p)(|\det {\rm D}_z \Psi_\delta(z)| - 1)| \le \varepsilon|\varphi(z)|^p\ ,
\ee
for $\delta > 0$ sufficiently small.  Incorporating \eqref{ineq_estim} into \eqref{ineq_phi}, we see that
$$
\int_{D_{0}}(v(\varphi(z))+2\varepsilon|\varphi(z)|^p)\,\md z\ge 0\ .
$$

It remains to find some $\beta=\beta(\varepsilon,\delta,\alpha) > 0$, such that for any admissible $\varphi$ in \eqref{half-ball}, the associated function $u=\delta^{-\frac{n}{p}}\varphi \circ \Psi^{-1}_\delta \circ \Phi_\delta$ is admissible as a test function in \eqref{ineq_u}, i.e., we need that 
$\|\mathcal{A}\varphi\|_{W^{-1,p}(D_{0};\R^d)} \leq \beta \|\varphi\|_{L^p(D_0; \R^m)}$ implies that $\|\mathcal{A}u\|_{W^{-1,p}(B(0,\delta) \cap \O;\R^d)} \le \alpha \|u\|_{L^p(B(0,\delta) \cap \O;\R^m)}$.

We calculate
\begin{align*}
&\norm{\cA\varphi}{W_0^{-1,p}(D_0;\R^d)}\\
 &= \sup_{\norm{w}{W_0^{1,p'}(D_0;\R^d)}
\leq 1} \sum_{i=1}^n\int_{D_0}A^{(i)}\varphi(z)\cdot \frac{\partial w(z)}{\partial z_i} \md x\\
&=\sup_{\|w\|\leq 1} \sum_{i=1}^n \int_{\Psi_\delta(D_0)}A^{(i)}\varphi(\Psi_\delta^{-1})\cdot\frac{\partial w}{\partial x'_i}(\Psi_\delta^{-1}(x'))\left|\det {\rm D} \Psi_\delta^{-1}(x')\right|\md x'\\
&=\begin{aligned}[t]
	\sup_{\|w\|\leq 1} \sum_{i=1}^n
	\int_{\frac{1}{\delta}(B(0,\delta)\cap 		\Omega)}\sum_{j=1}^d\left(A^{(i)}\varphi(\Psi_\delta^{-1}(x'))\right)_j 
	\Big({\rm D}\big(w(\Psi_\delta^{-1}(x'))\big)\cdot \big({\rm D}\Psi_\delta^{-1}(x')\big)^{-1}\Big)_{j,i}&\\
	\cdot\det |\rm D\Psi_\delta^{-1}(x')|&\md x'.
\end{aligned}
\end{align*}
Denoting $w_\delta:=w\circ \Psi_\delta^{-1}$, using the function $y$ as in \eqref{ineq64261} and the convergence of $\Psi_\delta^{-1}$ to the identity in $C^1(\overline{B(0,1)};\R^n)$,
we get
\begin{align*}
\norm{\cA\varphi}{W_0^{-1,p}(D_0;\R^d)}&{\geq} \frac{1}{2} \sup_{\norm{w_\delta}{W_0^{1,p'}(\Psi_\delta(D_0);\R^d)}
\leq 1} \sum_{i=1}^n\int_{\frac{1}{\delta}(B(0,\delta)\cap \Omega)}A^{(i)}y(x')\frac{\partial w_\delta(x')}{\partial x'_i}\md x'\\
&=\frac{1}{2}\sup_{\|w_\delta\|\leq 1} \sum_{i=1}^n\int_{B(0,\delta)\cap \Omega}A^{(i)}y(\delta^{-1}x)\frac{\partial w_\delta}{\partial x_i}(\delta^{-1}x)\md x\\
&=\frac{1}{2}\sup_{\|w_\delta\|\leq 1} \sum_{i=1}^n\int_{B(0,\delta)\cap\Omega}A^{(i)}\left(\delta^{n/p}u(x)\right)\delta\frac{\partial (w_\delta(\delta^{-1} x))}{\partial x_i} \md x
\end{align*}
for sufficiently small $\delta$. With $\eta_\delta(x):=\delta^{1-\frac{n}{p'}}w_{\delta}({\delta}^{-1}x)$ and due to
\[\norm{{\rm D}\eta_\delta}{L^{p'}(B(0,\delta)\cap\Omega;\R^d)}=\norm{{\rm D} w_\delta}{L^{p'}(\frac{1}{\delta}(B(0,\delta)\cap \Omega;\R^d)}\]
it follows that
\begin{align*}
\norm{\cA\varphi}{W_0^{-1,p}(D_0;\R^d)}&{\geq}\frac{1}{2}\sup_{\norm{\eta_\delta}{W_0^{1,p'}(B(0,\delta)\cap \Omega;\R^d)}\leq 1}\sum_{i=1}^n\int_{B(0,\delta)\cap \Omega}A^{(i)}u(x)\cdot \frac{\partial \eta_\delta(x)}{\partial x_i}\delta^{n}\md x\\
&=\frac{1}{2}\delta^{n}\norm{\cA u}{W^{-1,p}(B(0,\delta)\cap\Omega;\R^d)}.
\end{align*}

By a similar procedure as above, we compute
\begin{align*}
&\|u\|^p_{L^p(B(0,\delta) \cap \O; \R^m)} = \int_{B(0,\delta) \cap \O} |u(x)|^p \,\md x \\
& = \int_{\delta^{-1}(B(0,\delta) \cap \O)} |u(\Phi^{-1}_\delta(x')|^p |\det {\rm D}_{x'} \Phi^{-1}_\delta(x')| \,\md x' = \int_{\delta^{-1}(B(0,\delta) \cap \O)} |y(x')|^p \,\md x' \\
& = \int_{D_0} |y(\Psi_\delta(z))|^p |\det {\rm D}_z \Psi_\delta(z)| \,\md z \ge \frac{1}{2} \int_{D_0} |\varphi(z)|^p \,\md z =\frac{1}{2} \|\varphi\|^p_{L^p(D_0; \R^m)}\ .
\end{align*}
Hence, due to the assumption that $u$ is $\mathcal{A}$-qcb at $0$, we see that ($C_i>0$, $i=1,2,3$ are some constants) 
\[
\|\mathcal{A}\varphi\|_{W^{-1,p}(D_{0};\R^d)} \le C_1 \|\mathcal{A}u\|_{W^{-1,p}(B(0,\delta) \cap \O;\R^d)} \le C_2 \|u\|_{L^p(B(0,\delta) \cap \O; \R^m)} \le C_3 \|\varphi\|_{L^p(D_0; \R^m)}.
\]

``if'': The sufficiency of \eqref{half-ball} for $v$ to be $\mathcal{A}$-qcb at $0$ can be shown by analogous computations, instead of the (uniform) convergence of $\Psi_\delta$ one uses the (uniform) convergence of $\Psi_\delta^{-1}$ as $\delta \to 0$.
\hfill
$\Box$

\bigskip
Following the proof of Proposition~\ref{prop:aAqcbD}, we are also able to give an equivalent variant of $\cA$-qcb in the limit as $\delta\to 0$. 

\begin{proposition}\label{prop:xAqcbD}
Assume that $\O\subset\R^n$ has a boundary of class $C^1$ in a neighborhood of $x_0\in\partial\O$. Let $\nu_{x_0}$ be the outer unit normal to $\partial\O$ at $x_0$ and
\[
	D_{x_0}:=\{x\in B(0,1)\mid \,x\cdot\nu_{x_0}< 0\}.
\]
Then  $v\in \Hp(\R^m)$ is $\mathcal{A}$-qcb at $x_0$ if and only if
\be\label{half-ballB}
\begin{aligned}
	&\text{for every $\varepsilon>0$ there exists $\beta>0$ such that}\\
 	&\qquad\int_{D_{x_0}} v(\varphi(x))+\varepsilon|\varphi(x)|^p\,\md x\ge 0\   \\
	&\text{for every $\varphi\in L^p_0(B(0,\tfrac{1}{2});\R^m)$ with
	$\|\mathcal{A}\varphi\|_{W^{-1,p}(B(0,1);\R^d)}\le\beta \|\varphi\|_{L^p(D_{x_0};\R^m)}$.}
\end{aligned}
\ee
\end{proposition}

Unlike for strong $\cA$-qcb, it is possible to derive another version of $\cA$-qcb with periodic, precisely $\cA$-free test functions and a much more obvious relationship to $\cA$-quasiconvexity. Note however that instead of admitting test functions that are only ``almost'' $\cA$-free, we are then forced to work with a class that only ``almost'' has compact support (since $\gamma$ can be chosen arbitrarily small in \eqref{a-qcb-per} below).

%It again illustrates the implicit dependence on $\cA$-free extension: for the Cauchy--Riemann system, \eqref{a-qcb-per} below would be trivial, because all periodic and thus bounded holomorphic functions on $\CC$ are constant, and since $\gamma$ can be chosen small enough so that $|Q\setminus \frac{1}{2}Q|>\gamma^{\frac{1}{p}}|Q|$, $\varphi=0$ becomes the only admissible test function.

\begin{proposition}\label{prop:aqcb}\jknote{R11: Is it really superior? The test functions are $\mathcal{A}$-free, but they no longer have compact support. Compact support implies periodicity?}
Let $x_0\in\partial\O$, assume that $\partial \O$ is of class $C^1$ in a neighborhood of $x_0$, and define $Q=Q(x_0):=\{y\in \R^n\mid |y\cdot e_j|<1~\text{for}~j=1,\ldots,n\}$ and $Q^{-}:=\{y\in Q\mid y\cdot e_1<0\}$, where $e_1,\ldots,e_n$ of $\R^n$ is an orthonormal basis of $\R^n$ such that $e_1=\nu_{x_0}$, the unit outer normal to $\partial \O$ at $x_0$. Then $v \in \Hp(\R^m)$ is $\mathcal{A}$-qcb at $x_0$ if and only if
\be\label{a-qcb-per} 
\begin{aligned}
	&\text{for every $\varepsilon>0$, there exists $\gamma>0$ such that}\\
	&\qquad\int_{Q^-}  v(\varphi(x))+\varepsilon |\varphi(x)|^p \,\md x\ge 0\\
	&\text{for every $\varphi\in L^p_\#(Q;\R^m)$ with $\mathcal{A}\varphi=0$
	and $\|\varphi\|_{L^p(Q\setminus \frac{1}{2}Q;\R^m)}\leq \gamma \|\varphi\|_{L^p(Q;\R^m)}$.}
\end{aligned}
\ee
\end{proposition}
{\it Proof.}
``if'': We claim that \eqref{a-qcb-per} implies \eqref{half-ballB}. By $p$-homogeneity, it suffices to show the integral inequality in \eqref{half-ballB} for every $\varphi\in L^p_0(B(0,\tfrac{1}{2});\R^m)$ with $\|\varphi\|_{L^p}=1$ and $\|\cA \varphi\|_{W^{-1,p}}\leq \beta$, where $\beta=\beta(\varepsilon)$ is yet to be chosen. Below, the average of $\varphi$ is denoted by
\[
	a_\varphi := \frac{1}{|Q|}\int_Q \varphi(x)\,\md x.
\]
By Lemma~\ref{T} and Remark~\ref{rem:T}, $\|\varphi-a_\varphi-\TT \varphi\|_{L^p(Q;\R^m)}$ becomes arbitrarily small, provided that $\|\cA \varphi\|_{W^{-1,p}}\leq \beta$\jkerror{The domain seems to be critical, why is it missing} is small enough. In view of Lemma~\ref{lem:ucont} (uniform continuity of $u\mapsto v(u)$ and $u\mapsto |u|^p$, $L^p\to L^1$, on bounded sets in $L^p$), this means that for every $\varepsilon>0$, there exists a $\beta>0$ such that
$$
	\int_{Q^-} v(\varphi(x))+\varepsilon|\varphi(x)|^p\,\md x \geq \int_{Q^-} v(a_\varphi+\TT\varphi(x))+\frac{\varepsilon}{2}|a_\varphi+\TT\varphi(x)|^p\,\md x,
$$
and due to the inequality in \eqref{a-qcb-per} with $a_\varphi+\TT \varphi$ instead of $\varphi$, the right-hand side above is non-negative. Hence,
$$
	\int_{D_{x_0}} v(\varphi(x))+\varepsilon|\varphi(x)|^p\,\md x
	=\int_{Q^-} v(\varphi(x))+\varepsilon|\varphi(x)|^p\,\md x\geq 0.
$$

``only if'': Suppose that \eqref{half-ballB} holds. Let $\varepsilon>0$, and let $\varphi$ denote an admissible test function for \eqref{a-qcb-per}, i.e., $\varphi\in L^p_\#(Q_{x_0};\R^m)$ with $\mathcal{A}\varphi=0$ 	and $\|\varphi\|_{L^p(Q\setminus \frac{1}{2}Q;\R^m)}\leq \gamma \|\varphi\|_{L^p(Q;\R^m)}$, with some $\gamma$ still to be chosen. We may also assume that $\|\varphi\|_{L^p(Q)}=1$. Let $\eta\in C_0^\infty(Q;[0,1])$ be a fixed function such that $\eta=1$ on $\frac{1}{2}Q$ and $\eta=0$ on $Q\setminus \frac{3}{4}Q$. Observe that $	\|\varphi-\eta\varphi\|_{L^p(Q;\R^m)}\leq 2 \|\varphi\|_{L^p(Q\setminus \frac{1}{2}Q;\R^m)} \leq 2\gamma \|\varphi\|_{L^p(Q;\R^m)}$, whence
$$
	\|\varphi-\eta\varphi\|_{L^p(Q;\R^m)}\leq 2\gamma \|\varphi\|_{L^p(Q;\R^m)} 
	\leq \frac{2\gamma}{1-2\gamma} \|\eta \varphi\|_{L^p(Q;\R^m)} 
$$
In addition, there is a constant $C\geq 0$ depending on $\eta$ and $\cA$ such that
$$
	\|\cA(\eta\varphi)\|_{W^{-1,p}(Q;\R^d)}\leq C \|\varphi\|_{L^p(\frac{3}{4}Q\setminus \frac{1}{2}Q;\R^m)}
	\leq C\gamma \|\varphi\|_{L^p(Q;\R^m)} \leq \frac{C\gamma}{1-2\gamma} \|\eta\varphi\|_{L^p(Q;\R^m)} .
$$
Hence, for $\gamma$ sufficiently small, $\eta\varphi$ is an admissible test function for \eqref{half-ballB} (which we apply with $\varepsilon/2$ instead of $\varepsilon$), up to the fact that the support of $\eta\varphi$, which is contained in $\frac{3}{4}Q$, might be larger than $B(0,\frac{1}{2})$. This, however, can be easily corrected by a change of variables, rescaling by a fixed factor.
Consequently, 
$$
	\qquad\int_{Q_{x_0}^-}  v(\eta(x)\varphi(x))+\frac{\varepsilon}{2} |\eta(x)\varphi(x)|^p \,\md x\ge 0,
$$
and due to the uniform continuity shown in Lemma~\ref{lem:ucont}, we conclude that for $\gamma$ small enough,
$$
	\qquad\int_{Q_{x_0}^-}  v(\varphi(x))+\varepsilon |\varphi(x)|^p \,\md x\ge 0.
$$
\hfill
$\Box$

%%%%%%%%%%%%%%%%%%%%%%%%%%%%%%%%%%%%%%%%%%%%%%%%%%%%%%%%%%%%%%%%%%%%%%%%%%%%%%%%%%%%%%%%%%%%%%

We now focus on the link between (strong) $\cA$-quasiconvexity at the boundary and weak lower semicontinuity along (asymptotically) $\cA$-free sequences. 

\section{Link to weak lower semicontinuity}
\subsection{Asymptotically $\cA$-free sequences}

\begin{proposition}\label{prop:awlsc}
Let $h_\infty\in C(\bar\O;\Hp(\R^m))$. Then $I_\infty(u):=\int_\O h_\infty(x,u(x))\,\md x$ is weakly sequentially lower semicontinuous  along asymptotically $\cA$-free sequences in $L^p(\O;\R^m)$ if and only if 

\noindent (i) $h_\infty$ is strongly-$\mathcal{A}$-qcb at every $x_0\in\partial\O$ and 

\noindent (ii) $h_\infty(x,\cdot)$ is $\mathcal{A}$-quasiconvex at almost every $x\in\O$.
\end{proposition} 

{\it Proof.}
``only if'': We show that  strongly-$\mathcal{A}$-qcb at $x_0\in\partial\O$ is a necessary condition; the necessity of (ii) is well known. Suppose that $h_\infty$ is not  strongly-$\mathcal{A}$-qcb at $x_0\in\partial\O$. This means that there is $\varepsilon>0$ such that for every $k\in\N$ there exists $u_k\in L^p_0(B(x_0,\frac{1}{k});\R^m)$ with  $\|\mathcal{A}u\|_{W^{-1,p}(\O;\R^d)}\le \frac{1}{k} \|u_k\|_{L^p(\O;\R^m)}$ 
and 
\[
	\int_{B(x_0, \frac{1}{k})\cap\O} h_\infty(x,u_k(x))+\varepsilon|u_k(x)|^p\,\md x<0\ .
\]
In particular, $u_k$ cannot be the zero function.
Denote 
\[
	\hat u_k:=	u_k/\|u_k\|_{L^p(B(x_0,\frac{1}{k})\cap \O;\R^m)} = u_k/\|u_k\|_{L^p(\O;\R^m)}.
\]
Then $\hat u_k\in L^p_0(\O;\R^m)$ with $\|\hat u_k\|_{L^p}= 1 $ and $\|\mathcal{A}\hat u_k\|_{W^{-1,p}(\O;\R^d)} \leq 1/k$. In addition, $\hat u_k$ vanishes outside of $B(x_0,\frac{1}{k})$, so that $\hat u_k\to 0$ in measure and weakly in $L^p(B(x_0,1);\R^m)$. However, 
$$ 
\liminf_{k\to\infty} \int_{\O}h_\infty(x,\hat u_k(x))\,\md x \leq-\varepsilon<0=\int_{\O} h_\infty(x,0)\,dx\ .
$$
This means that $u\mapsto\int_\O h_\infty(x,u(x))\,\md x$ is not lower semicontinuous along $\{\hat{u}_k\}$.

``if'': Let us now prove the sufficiency. Let $\{u_k\}_{k\in\N}\subset L^p(\O;\R^m)$ be an asymptotically $\cA$-free sequence weakly converging to some $u$ in $L^p$. As a first step, we assume that in addition, $\{u_k\}$ is purely concentrating in the sense that $u_k\rightharpoonup 0$ in $L^p(\O;\R^m)$ and  $\mathcal{L}^n(\{x\in\O;\,u_k(x)\ne 0\})\to 0$ as $k\to\infty$. It suffices to show that every subsequence of $\{u_k\}$ admits another subsequence along which $I$ is lower semicontinuous. Using DiPerna-Majda measures as in \eqref{special} in the Appendix, and  we get that for every $\delta>0$, up to a subsequence,
\be \label{DPMlimit}
\begin{aligned}
& \lim_{k\to\infty} \int_{B(x_0,\delta)\cap\O} h_\infty(x,u_k(x))\,\md x\\
&\qquad = \int_{\overline{B(x_0,\delta)\cap\O}}\int_{\b_{\cal S}\R^m\setminus\R^m} \frac{h_\infty(x,s)}{1+|s|^p}\md\lambda_x(s)\md\pi(x)\ 
\end{aligned}
\ee
for some $(\pi,\lambda)\in\cdm$. 

In the following, we only consider those $\delta>0$ for which $\pi(\partial B(x_0,\delta)\cap \bar\O)=0$, which is certainly true for a dense subset. Let $\{\eta_\ell\}_{\ell\in\N}\subset C_0^\infty(B(x_0,\delta))$ such that $0\le\eta_\ell\le 1$ and $\eta_\ell\to\chi_{B(x_0,\delta)}$  as $\ell\to\infty$. Here,  $\chi_{B(x_0,\delta)}$ is the characteristic function of $B(x_0,\delta)$ in $\R^n$ and $x_0\in\partial\O$. By Lemma~\ref{lem:truncate}, $\mathcal{A}(\eta_\ell u_k)\to 0$ in $W^{-1,p}(\O;\R^d)$ as $k\to\infty$, for fixed $\ell$. Take $\varepsilon>0$, $x_0\in\partial\O$, $\alpha,\delta>0$ as in Definition~\ref{def:aaqcb} and set $w_k:=\eta_{\ell(k)}u_k$, where $\ell(k)$ tends to $\infty$ sufficiently slowly as $k\to\infty$ so that $\cA w_k\to 0$ in $W^{-1,p}(\O;\R^d)$ and reasoning as in \cite[Appendix]{ifmk}, using that $\pi(\partial B(x_0,\delta)\cap \bar\O)=0$, we see that $\{w_k\}$ also generates $(\pi,\lambda)$, at least on $\overline{B(x_0,\delta)}\cap \Omega$. If $w_k$ strongly converges to zero in $L^p$,
\be\label{limalongwknonneg}
\begin{aligned}
	0&\le \lim_{k\to\infty}\int_{B(x_0,\delta)\cap\O} h_\infty(x,w_k(x))+\varepsilon|w_k(x)|^p\,\md x,
\end{aligned}
\ee
by continuity (in that case, we even get equality). Otherwise, a subsequence of $\{w_k\}$ (not relabeled) is bounded away from zero in $L^p$, and since $\cA w_k\to 0$ in $W^{-1,p}$, this implies that $\|\cA w_k\|_{W^{-1,p}}\leq \alpha \|w_k\|_{L^p}$, at least for $k$ large enough. Hence, $w_k$ is admissible as a test function in \eqref{aa-qcb}, and we end up again with \eqref{limalongwknonneg}. The right-hand side of \eqref{limalongwknonneg} can be expressed using \eqref{special}:
\[
\begin{aligned}
	&\lim_{k\to\infty}\int_{B(x_0,\delta)\cap\O} h_\infty(x,w_k(x))+\varepsilon|w_k(x)|^p\,\md x\\
 &\qquad = \int_{\overline{B(x_0,\delta)\cap\O}}\int_{\b_{\cal R}\R^m\setminus\R^m} \frac{h_\infty(x,s)+\varepsilon |s|^p}{1+|s|^p}\md\lambda_x(s)\md\pi(x)\ .
\end{aligned}
\]
Hence,
$$
0\le \pi(\overline{B(x_0,\delta)\cap\O})^{-1}\int_{\overline{B(x_0,\delta)\cap\O}}\int_{\b_{\cal R}\R^m\setminus\R^m} \frac{h_\infty(x,s)+\varepsilon |s|^p}{1+|s|^p}\md\lambda_x(s)\md\pi(x)\ .$$
Therefore, by the  Lebesgue-Besicovitch differentiation theorem (see \cite{evans}, e.g.) and by taking into account that $\varepsilon>0$ is arbitrary we get that for $\pi$-almost every $x_0\in\partial\O$
$$
0\le\int_{\b_{\cal R}\R^m\setminus\R^m} \frac{h_\infty(x_0,s)}{1+|s|^p}\md\lambda_{x_0}(s)\ .
$$
This together with Theorem~\ref{th-app} and \eqref{special} implies that the inner integral on the right-hand side of  \eqref{DPMlimit} is nonnegative
for $\pi$-almost every $x_0\in\bar\O$. As a consequence, $I_\infty$ is lower semicontinuous along $\{u_k\}$, i.e., all purely concentrating sequences. By Theorem~\ref{thm:wlscaAfree} and Remark~\ref{replacement} (ii), we conclude that $u\mapsto \int_\O h(x,u(x))\,\md x$ is weakly lower semicontinuous along arbitrary asymptotically $\mathcal{A}$-free sequences.
\hfill
$\Box$ 

\bigskip

In view of Remark~\ref{replacement}, our results obtained so far can be summarized as follows.
\begin{theorem}\label{thm:wlscasympt} Suppose that  $\O\subset\R^n$ be a bounded domain with $\mathcal{L}^n(\partial\O)=0$.
Let $1<p<+\infty$, and let $h:\bar\O\times\R^m\to\R$ be continuous and such that $h(x,\cdot)\in\ups$ for all $x\in\bar\O$, with recession function $h_\infty\in C(\overline{\O};\Hp)$.
Then $I$ is weakly
lower semicontinuous along asymptotically $\cA$-free sequences if and only if 

\noindent (i) $h(x,\cdot)$ is $\mathcal{A}$-quasiconvex for almost all $x\in\O$;

\noindent (ii) $h_\infty$ is strongly $\mathcal{A}$-quasiconvex at the boundary for all $x_0\in\partial\O$.
\end{theorem}

\bigskip

\subsection{Genuinely $\cA$-free sequences}

We now focus on weak lower semicontinuity along sequences $\{u_k\}$ that satisfy $\cA u_k=0$ for each $k\in\N$. Since a substantial part of the arguments in this context is analogous to the ones in the preceding subsection, we do not always give full proofs. The main difference is that for the link to $\cA$-quasiconvexity at the boundary ($\cA$-qcb) as introduced in Definition~\ref{def:aqcb}, more precisely, for its sufficiency, we rely on an extension property:

\begin{definition}[$\cA$-free extension domain]\label{def:Aext}
We say that $\O$ is an \emph{$\cA$-free extension domain} if there exists a larger domain $\O'$ with $\O\subset\subset \O'$ and an associated \emph{$\cA$-free extension operator}, i.e., a bounded linear operator $E:L^p(\O;\R^m)\cap \ker\to L^p(\O';\R^m)\cap \ker$ such that $Eu=u$ on $\O$.
\end{definition}

As mentioned before, the existence of an $\cA$-free extension operator not only depends on the smoothness of $\partial\O$, but also on $\cA$ itself. On the one hand, if $\partial\O$ is Lipschitz,  extension operators are available for $\cA=\rm{curl}$ and $\cA=\rm{div}$ (essentially using a partition of unity and extension by a suitable reflection), but on the other hand, if we choose $\cA$ to be the differential operator of the Cauchy--Riemann system ($n=m=2$, identifying $\CC$ with $\R^2$), no such extension operator exists even for very smooth domains, since holomorphic functions with singularities at the boundary of $\O$ can never be extended to holomorphic functions on a larger set including the singular point\footnote{In terms of integrability, the weakest possible point singularity of an elsewhere holomorphic function locally behaves like $z\mapsto 1/z$ ($z\in\CC\setminus\{0\}$), which is not in $L^p(\O)$ if $p\geq 2$, $0\in \partial \O$ and $\partial \O$ is smooth in a neighborhood, but using an appropriately weighted series of singular terms, each with a singularity slightly outside $\O$, accumulating at a boundary point, examples in $L^p$ are possible for arbitrary $1\leq p<\infty$.}.

With the help of the extension property and the projection $\TT$  of Lemma~\ref{T},
Proposition~\ref{prop:awlsc} can be adapted to the setting of genuinely $\cA$-free sequences:
\begin{proposition}\label{prop:wlsc}\jknote{R8: It is 'iff', of course we can replace it.}
Suppose that $\O$ is an $\cA$-free extension domain and let $h_\infty\in C(\bar\O;\Hp(\R^m))$. Then $I_\infty(u):=\int_\O h_\infty(x,u(x))\,\md x$ is weakly sequentially lower semicontinuous along $\cA$-free sequences in $L^p(\O;\R^m)$ if and only if 

\noindent (i) $h_\infty$ is $\mathcal{A}$-qcb at every $x_0\in\partial\O$ and 

\noindent (ii) $h_\infty(x,\cdot)$ is $\mathcal{A}$-quasiconvex at almost every $x\in\O$.
\end{proposition} 
{\it Proof.} 
``only if'': Again, necessity of (ii) is well known. If $h_\infty$ is not $\mathcal{A}$-qcb at a point $x_0\in \partial \O$, as in the proof of Proposition~\ref{prop:awlsc} we obtain an $\varepsilon>0$ and a sequence $\{\hat{u}_k\} \subset L^p_0(B(x_0,\frac{1}{k});\R^m)$ with $\|\hat u_k\|_{L^p(\O;\R^m)}= 1$ such that 
$$ 
\liminf_{k\to\infty} \int_{\O}h_\infty(x,\hat u_k(x))\,\md x \leq-\varepsilon<0=\int_{\O} h_\infty(x,0) \,\md x,
$$
and $\|\mathcal{A}\hat u_k\|_{W^{-1,p}(\R^n;\R^d)} \leq 1/k$. Each $\hat{u}_k$ can be interpreted as a $Q$-periodic function $\hat{u}_k^\#$ with respect to a cube $Q$ compactly containing $\O\cup B(x_0,1)$, by first extending $\hat{u}_k$ by zero to the rest of $Q$ and then periodically to $\R^n$. 
We denote its cell average by
\[
	a_k:=\frac{1}{|Q|} \int_Q \hat{u}_k\,\md x.
\]
By Remark~\ref{rem:T}, we infer that $\|\mathcal{A}\hat u_k^\#\|_{W_\#^{-1,p}(\R^n;\R^d)} \leq C/k$ with a constant $C\geq 0$ independent of $k$.
The projection of Lemma~\ref{T} now yields the sequence $\{\TT\hat{u}_k^\#\} \subset L^p_\#(R^n;\R^m)\cap \ker$, which satisfies $\|a_k+\TT\hat{u}_k^\#-\hat{u}_k\|_{L^p(Q;\R^m)}\to 0$ as $k\to\infty$. Consequently, $a_k+\TT\hat{u}_k^\#\rightharpoonup 0$ weakly in $L^p$ just like $\hat{u}_k$, and due to Lemma~\ref{lem:ucont} (uniform continuity on bounded subsets of $L^p$), 
$$ 
\liminf_{k\to\infty} \int_{\O}h_\infty(x,a_k+\TT\hat u_k^\#(x))\,\md x \leq-\varepsilon<0=\int_{\O} h_\infty(x,0)\,\md x.
$$
Hence, $I_\infty$ is not lower semicontinuous along the $\cA$-free sequence $\{a_k+\TT\hat u_k^\#\}$.

``if'': The argument is completely analogous to that of Proposition~\ref{prop:awlsc}, using Theorem~\ref{wlsc1} instead of Theorem~\ref{thm:wlscaAfree}. Observe that due to the extension operator, any given sequence $\{u_k\}$ along which we want to show lower semicontinuity is defined and $\cA$-free on some set $\O'\supset\supset \O$. Hence, after the truncation argument of Proposition~\ref{prop:awlsc}, we now end up with an admissible test function for Definition~\ref{def:aqcb} (see also Remark~\ref{rem:Aqcb-test}).
\hfill
$\Box$

\bigskip

We arrive at the analogous main result:

\begin{theorem}\label{thm:wlscAfree}
Suppose that  $\O\subset\R^n$ be a bounded $\mathcal{A}$-free extension domain with $\mathcal{L}^n(\partial\O)=0$. 
Let $1<p<+\infty$, and let $h:\bar\O\times\R^m\to\R$ be continuous and such that $h(x,\cdot)\in\ups$ for all $x\in\bar\O$, with recession function $h_\infty\in C(\overline{\O};\Hp)$.
Then $I$ is sequentially weakly
lower semicontinuous along  $\cA$-free sequences if and only if 

\noindent (i) $h(x,\cdot)$ is $\mathcal{A}$-quasiconvex for almost all $x\in\O$;

\noindent (ii) $h_\infty$ is $\mathcal{A}$-quasiconvex at the boundary for all $x_0\in\partial\O$.
\end{theorem}
\begin{remark}
In general, the continuity of $h_\infty$ in $x$ cannot be dropped in Theorem~\ref{thm:wlscAfree}. For a counterexample in the gradient case ($\cA=$curl) see \cite[Section 4]{kroemer}.
\end{remark}

\section{Concluding remarks}\label{sec:concluding}
\subsection{$\cA$-free versus asymptotically $\cA$-free sequences}\label{subsec:approx}

Clearly, weak lower semicontinuity along asymptotically $\cA$-free sequences implies weak sequential lower semicontinuity for the functional restricted to $\ker$. We do not know whether or not the converse is true in general.\jknote{4R1: This should answer his question, we do not know of a counterexample, otherwise the converse would not be true in general} However, it holds at least in some special cases. More precisely, it suffices to have an extension property in the following sense. It trivially implies the $\cA$-free extension property mentioned in Definition~\ref{def:Aext} (but the converse is not clear there, either): 

\begin{definition}\label{def:afe}[asymptotically $\cA$-free extensions]\jknote{Fixed spelling mistake}
We say that $\O$ has the $\cA$-($L^p$,$W^{-1,p}$) extension property if there exists a domain $\Lambda$ with $\bar\O\subset\subset \Lambda$ such that
for every $u\in L^p(\Omega;\R^m)$, there is an extension $v\in L^p(\Lambda;\R^m)$ of $u$ which satisfies
$$
	\|v\|_{L^p(\Lambda;\R^m)} \leq C \|u\|_{L^p(\O;\R^m)}\quad\text{and}\quad
	\|\cA v\|_{W^{-1,p}(\Lambda;\R^d)} \leq C \|\cA u\|_{W^{-1,p}(\O;\R^d)},
$$
where $C\geq 0$ is a suitable constant only depending on $\Lambda$, $\O$, $p$ and $\cA$.
\end{definition}

If this holds, we can always reduce asymptotically $\cA$-free sequences to genuinely $\cA$-free sequences with arbitrarily small error in $L^p$. The argument can be sketched as follows: For a given approximately $\cA$-free sequence $u_k\rightharpoonup u$ along which we want to show lower semicontinuity, it is possible to truncate the extension of $u_k-u$, multiplying with a cut-off function which is $1$ on $\Omega$ and makes a transition down to zero in $\Lambda\setminus \Omega$ (this cannot be done inside, because $u_k$ might concentrate a lot of mass near the boundary, and cutting off inside could then significantly alter the limit of the functional along the sequence). The modified sequence is still asymptotically $\cA$-free due to Lemma~\ref{lem:truncate}, and since it is compactly supported in $\Lambda$ by construction, we can further extend it periodically to $\R^n$, with a sufficiently large fundamental cell of periodicity containing the support of the cut-off function. We thus end up in the periodic setting where we can project onto $\cA$-free fields with controllable error, using Lemma~\ref{T}. 

Clearly, the $\cA$-($L^p$,$W^{-1,p}$) extension property implies the standard $\mathcal{A}$-free extension property introduced in Def~\ref{def:Aext}, and if the former holds, then $\cA$-quasiconvexity at the boundary and strong $\cA$-quasiconvexity at the boundary are equivalent.
Even for smooth domains, the $\cA$-($L^p$,$W^{-1,p}$) extension property depends on $\cA$ (and possibly on $p$). 
For instance, it holds for $\cA={\rm div}$ on domains of class $C^1$ using local maps and extension by an appropriate reflection for flat pieces of the boundary, but not for all $\cA$. In particular, it fails for the Cauchy-Riemann system, just like the weaker $\mathcal{A}$-free extension property introduced in Def~\ref{def:Aext}. 
Interestingly, the case $\cA=$curl seems to be nontrivial: the $\cA$-($L^p$,$W^{-1,p}$) extension property for $\cA=$curl does hold for $n=2$ (the $2d$-curl and the $2d$-divergence are the same operators up to a fixed rotation), but if $n\geq 3$, we do not know. For a flat piece of the boundary, the natural extension for almost curl-free fields would of course also be by reflection, i.e., the one corresponding to an even extension of the scalar potential across the boundary (even in direction of the normal), but in this case,
the required estimate in $W^{-1,p}$ for the curl seems to be nontrivial, if true at all. The problem appears for those of components of the curl that only contain partial derivatives in tangential directions, precisely the ones that ``naturally'' get extended to even functions, say, $\partial_2 u_3-\partial_3 u_2$, if the normal to the boundary (locally) is the first unit vector.

The situation for smooth domains is summarized in the table below: 
\begin{center}
	\begin{tabular}{|l|l|l|}
		\hline
		 $\cA$ & strong $\cA$-qcb ~$\Leftrightarrow$~$\cA$-qcb & Extension property of Def.~\ref{def:afe}  \\
		\hline
		div ($n\in\N$) & true & true\\
		Cauchy-Riemann ($n=2$) & false & false\\ 
		curl ($n=2$) & true  & true  \\
		curl ($n>2$) & open & open\\
		\hline
	\end{tabular}
\end{center}
Although the second and the third column in the table coincide we do not know whether the existence of the extension in the sense of Def.~\ref{def:afe} is really equivalent to  $\cA$-qcb~$\Rightarrow$~strong $\cA$-qcb. In view of the constant rank condition \eqref{rank} which makes it hard to characterize the class of the admissible operators $\cA$ beyond a few examples, a systematic analysis for all $\cA$ seems to be out of reach. 

%even if the piece of the boundary of a domain in $\CC\cong \R^2$ across which we want to extend is already flat: 

\subsection{The gradient case and classical quasiconvexity at the boundary}\label{subsec:gradient}

If $\varphi\in\ker$ then \eqref{half-ball} as well as \eqref{half-ballB} implies that  $\int_{D_{x_0}} v(\varphi(x))\,\md x\ge 0$. 
For $\cA=\operatorname{curl}$, the differential constraint can also be encoded using potentials: If $\varphi\in L^p$ and $\operatorname{curl} \varphi=0$ on the simply connected domain $D_{x_0}$, then there exists a potential vector field $\Phi\in W^{1,p}$ with $\varphi=\nabla \Phi$, and if $\varphi=0$ on $D_{x_0}\setminus B(0,\frac{1}{2})$, then $\Phi$ inherits this property up to an appropriate choice of the constants of integration. Hence, we get that
\be\label{half-ballgrad}
\begin{aligned}
 	\int_{D_{x_0}} v(\nabla \Phi(x))\,\md x\ge 0\ \qquad	&\text{for every $\Phi\in W_0^{1,p}(B(0,\tfrac{1}{2});\R^m)$.}
\end{aligned}
\ee
Taking into account that for $p$-homogeneous $v$, $v(0)=0$ and ${\rm D}v(0)=0$, the latter condition is the so-called quasiconvexity at the boundary \cite{bama} (at the zero matrix). 

The converse, that is, going back from \eqref{half-ballgrad} to either \eqref{half-ball} or \eqref{half-ballB}, is not so obvious, however. Nevertheless, in  case of \eqref{half-ballB}, this is true as a consequence of known characterizations of weak lower semicontinuity, on the one hand our Proposition~\ref{prop:wlsc} and the other hand Theorem 1.6 in \cite{kroemer}. (A proof directly working with the two conditions is also possible, although slightly more technical.)
%Both results apply to functionals of the form $U\mapsto \int_{\Omega_{x_0}} \eta(x) v(\nabla U(x))\,dx$, where $v$ is continuous and $p$-homogeneous, $\eta\in C^\infty(\R^n;[0,1])$ is supported in a small neighborhood of $x_0\in \partial\O_{x_0}$, and $\O_{x_0}$ is a $C^1$-domain chosen in such a way that  in a neighborhood of $x_0$ containing the support of $\eta$ on $\partial \O_{x_0}$, $\partial O_{x_0}$ is a hyperplane with constant normal $\nu_{x_0}$. (A proof directly working with the two conditions is also possible, although slightly more technical.)
%As to the question whether \eqref{half-ballgrad} implies \eqref{half-ball} for $\cA=\operatorname{curl}$, we suspect that at least for $n\geq 3$, this is not true in general, but our attempts of constructing an example of a function $v$ proving this so far did not succeed.

\subsection{Examples for the case of higher order derivatives}\label{subsec:hessian}

The following example shows that $I(u):=\int_\O {\rm det}\,\nabla^2 u(x)\,\md x$ is not weakly lower semicontinuous on $W^{2,2}(\O)$. 
Consequently, the determinant is not $\mathcal{A}$-qcb for suitably defined $\mathcal{A}$. As to the definition of $\cA$, we recall \cite{fomu}: The functional $I$ fits into our framework, if instead of $\nabla^2 u$, we define $I$ on fields $v=(v)_{ij}$, $1\leq i\leq j\leq n$, in $L^2$, satisfying $\cA v:=\operatorname{curl}~v=0$, with the understanding that for each $x$, $v(x)$ (the upper triangular part of a matrix) is identified with a symmetric matrix in $\R^{n\times n}$ still denoted $v$, both for the application of the (row-wise) $\operatorname{curl}$ and the evaluation of $I$, where $\nabla^2 u$ is replaced by $v$. One can check that $\cA v=0$ if and only if there exists a scalar-valued $u\in W^{2,2}$ with $v=\nabla^2 u$, at least as long as the domain is simply connected.

\begin{example}
  Consider $\O:=(-1,1)^2$ and for $F\in\R^{2\times 2}$ the function  $v_\infty(F):={\rm det}\, F$ and  the operator $\mathcal{A}$ such that $\mathcal{A}w=0$ if and only if for some $u\in W^{2,2}(\O)$, $w$ is the upper (or lower) triangular part of $\nabla ^2u$, which takes values in the symmetric matrices; cf.~\cite[Example 3.10(d)]{fomu}. Here $\nabla^2u$ denotes the Hessian matrix of $u$. Then $v_\infty$ is not $\mathcal{A}$-qcb. Indeed, take $u\in W_0^{2,2}(\O)$ extended by zero to the whole $\R^2$.  Define $u_k(x):= k^{-1}u(kx)$. Then $u_k\rightharpoonup 0$ in $W^{2,2}(\O)$. We have that 
\be
\lim_{k\to\infty}\int_{(0,1)\times(-1,1)}{\rm det}\, \nabla^2 u_k(x)\,\md x=  \int_{(0,1)\times(-1,1)}{\rm det}\, \nabla^2 u(y)\,\md y\ .
\ee 
Hence, it remains to find $u$ for which the integral on the right-hand side is negative which is certainly possible.
%Let $u(x_1,x_2):=f(x_1)g(x_2)$ where $f,g:[-1,1]\to\R$ are smooth and such that $g(\pm 1)=g'(\pm 1)=f(1)=f'(1)=0$, $f'(0)f(0)>0$, and $g'$ does vanish identically.
%Then 
%\[
%\begin{aligned}
%\int_{(0,1)\times(-1,1)}{\rm det}\, \nabla^2 u(y)\,\md y
%&= \int_{(0,1)\times(-1,1)}f(x_1)g(x_2)f''(x_1)g''(x_2)-f'(x_1)^2g'(x_2)^2\,\md x\\
%&=\begin{aligned}[t] 
	%&\int_{(0,1)\times(-1,1)}f'(x_1)^2g'(x_2)^2\,\md x+[f'(x)f(x)]^1_0\int_{-1}^1 g''(x_2)g(x_2)\,\md x_2\\
	%&-[g'(x)g(x)]^1_{-1}\int_{0}^1 f'(x)^2\,\md x_1-\int_{(0,1)\times(-1,1)}f'(x_1)^2g'(x_2)^2\,\md x
%\end{aligned}\\
%&=-f'(0)f(0)\int_{-1}^{1}g'(x_2)^2\,\md x<0\ .
%\end{aligned}
%\]
\end{example}

In the next example, we isolate a function which is $\mathcal{A}$-quasiconvex at the boundary.

\begin{example}
Consider $\O:=B(0,1) \subset \R^3$ and $\mathcal{A}$ such that $\mathcal{A}w = 0$ if and only if $w = \nabla^2 u$ for some $u \in W^{2,2}(\O)$, and the mapping $h(x,F):=a(x)\cdot({\rm Cof} F)\nu(x)$, where $a\in C(\bar\O;\R^3)$ is arbitrary and $\nu(x)\in C(\bar\O)$ coincides with the outer unit normal to $\partial\O$ for $x\in\partial\O$. Notice that by definition of the Cofactor matrix ($({\rm Cof} F)_{ij}$ is $(-1)^{i+j}$ times the determinant of the $2\times 2$ submatrix of $F$ obtained by erasing the $i$-th row and $j$-th column), $({\rm Cof} \nabla u(x))\nu(x)$ effectively only depends on directional derivatives of $u$ in directions perpendicular to $\nu(x)$.

For this $h$,
$$ 
\int_\O h(x,\nabla^2u_k(x))\,\md x\to\int_\O h(x,\nabla^2u_0(x))\,\md x
$$
whenever $u_k\rightharpoonup u_0$ in $W^{2,2}(\O)$. 

To see that consider $z_k:=\nabla u_k$ for $k\in\N\cup\{0\}$. Then $\{z_k\}\in W^{1,2}(\O;\R^3)$ and the result follows from\cite{kkk}.
\end{example}

\appendix
\section{Appendix}
\subsection{DiPerna-Majda measures}
In what follows, $\rca(\bar\O)$ denotes the space of Radon measures on $\bar\O$.
Consider the following complete (i.e. containing constants, separating points
from closed subsets and closed with respect to the supremum norm), 
separable (i.e. containing a dense countable subset) ring $\mathcal{S}$ of
continuous bounded functions from  $\R^{m}$ into $\R$ defined as 
\be\label{spherecomp}
\begin{aligned}
	{\cal S}:=\bigg\{ v_0\in C(\R^m) \,\bigg|\, \mbox{ there exist } c\in\R\ ,\  v_{0,0}\in C_0(\R^m),\mbox{ and } v_{0,1}\in C(S^{m-1}) \mbox{ s.t. } &\\
 v_0(s) = c+ v_{0,0}(s)+v_{0,1}\left(\frac{s}{|s|}\right)
\frac{|s|^p}{1+|s|^p}\mbox { if $s\ne 0$ and }  v_0(0)=c+v_{0,0}(0)&\bigg\} ,
\end{aligned}
\ee
where $S^{m-1}$ denotes the $(m-1)$-dimensional unit sphere in $\R^m$. 
 It is known that there is a
one-to-one correspondence ${\cal R}\mapsto\b_{\cal R}\R^{m}$ between such
rings and a (metrizable) compactification $\b_{\cal R}\R^{m}$ of $\R^{m}$ \cite{engelking}; for ${\cal R}={\cal
S}$, $\b_{\cal S}\R^{m}$ is obtained by adding a sphere to $\R^{m}$ at infinity. More precisely, $\b_{\cal
S}\R^m$ is homeomorphic to the closed unit ball $\overline{B(0,1)}\subset \R^m$ via the mapping $f:\R^m\to B(0,1)$, $f(s):=s/(1+|s|)$ for all $s\in\R^m$. Note that $f(\R^m)$ is dense in $\overline{B(0,1)}$.
For simplicity, we
will not distinguish between $\R^{m}$ and its image in $\b_{\cal S}\R^{m}$.

DiPerna and Majda \cite{diperna-majda} proved the following theorem:  

\begin{theorem}
Let $\O$ be an open 
domain in $\Rn$ with $\mathcal{L}^n(\partial\O)=0$,  and 
 let $\{y_k\}_{k\in\N}\subset L^p(\O;\R^{m})$, with $1\le p<+\infty$, be bounded.
 Then there exists a subsequence (not relabeled), a positive Radon measure $\pi\in\rca(\bar\O)$ and a family of probability measures on $\b_{\cal S}\R^{m}$ $\lambda:=\{\lambda_x\}_{x\in\bar\O}$  such that 
for all $h_0\in C(\bar\O\times\b_{\cal S}\R^{m})$ it holds that  
\be\label{basic}
\lim_{k\to\infty}\int_\O h_0(x,y_k(x))(1+|y_k(x)|^p)\d x\ =
\int_{\bar\O}\int_{\b_{\cal S}\R^{m}}h_0(x,s)\d\lambda_x(s)\d\s(x)\ . 
\ee
\end{theorem}

\bigskip

If \eqref{basic} holds we say that $\{y_k\}$ generates $(\pi,\lambda)$ and we denote the set of all  such pairs of measures generated by some sequence in $L^p(\O;\R^m)$ by $\cdm$.

For any $h(x,s):=h_0(x,s)(1+|s|^p)$ with $h_0\in C(\bar\O\times \b_{\cal S}\R^{m})$ then  there exists a continuous and positively $p$-homogeneous function $h_\infty:\bar\O\times\R^m\to\R$, i.e., $h_\infty(x,t s)=t^p h_\infty(x,s)$ for all $t\ge 0 $, all $x\in\bar\O$, and  $s\in\R^m$, such that 
\be\label{recession}
\lim_{|s|\to\infty}\frac{h(x,s)-h_\infty(x,s)}{|s|^p}=0\ .
\ee

 It is already mentioned in \cite{ifmk,kruzik} that if $\{y_k\}\subset L^p(\O;\R^m)$ is bounded and $\mathcal{L}^n(\{x\in\O;\,y_k(x)\ne 0\})\to 0$ as $k\to\infty$ then \eqref{basic} can be replaced by 
 \be\label{special}
 \lim_{k\to\infty}\int_\O h_\infty(x,y_k(x))\d x\ =
\int_{\bar\O}\int_{\b_{\cal S}\R^{m}\setminus\R^m}\frac{h_\infty(x,s)}{1+|s|^p)}\d\lambda_x(s)\d\s(x)\ ,
 \ee 
where $(x,s)\mapsto h_0(x,s):=h_\infty(x,s)/(1+|s|^p)$ belongs to $C(\bar\O\times \b_{\cal S}\R^{m})$.  

\bigskip

The following theorem is a direct consequence of \cite[Thms.~2.1, 2.2]{ifmk}.

\begin{theorem}\label{th-app}
 Let $\{y_k\}\subset L^p(\O;\R^m)\cap\ker$ generates $(\s,\lambda)\in\cdm$  and let $y_k\to 0$ in measure.  Then 
 for $\s$-almost
every $x\in\O$ and all $h\in C(\bar\O;\Hp(\R^m))$ such that $h(x,\cdot)$ is $\mathcal{A}$-quasiconvex for all $x\in\bar\O$  it holds
that \be\label{rem1}
 0\le  \int_{\beta_{{\cal S}}\R^{m}\setminus\R^{m}}\frac{h(x,s)}{1+|s|^p}\md\lambda_x(s)\ .
    \ee

\end{theorem}

\subsection{Uniform continuity properties of the functional}
The following lemma essentially allows us to modify sequences inside $I$ as long as the 
modified sequences approaches the original one in the norm of $L^p$.

%Given an $\cA$-free extension operator, especially if it maps to $L_{\#}^p(\R^n,\R^m)\cap\ker$, the projection $\TT$ of Lemma~\ref{T} 
%can be used to regain $\cA$-free fields after modifications leading outside $\ker$, for instance multiplication with cut-off functions which frequently appears in the proofs. To make sure that such technical operations do not affect the value of the functional too much, we rely on its uniform continuity on bounded subsets of $L^p$:

\begin{lemma}\label{lem:ucont}
Let $h_\infty\in C(\bar\O;\Hp(\R^m))$. Then for any pair $\{u_k\}$, $\{v_k\}$ of bounded sequences in $L^p(\O;\R^m)$ such that $u_k-v_k\to 0$ strongly in $L^p$, $h_\infty(\cdot,u_k(\cdot))-h_\infty(\cdot,v_k(\cdot))\to 0$ strongly in $L^1$.
\end{lemma}
{\it Proof.} 
For $\delta>0$ let 
\[
	A_k(\delta):=\{x\in\O:|u_k(x)-v_k(x)|\geq \delta (|u_k(x)|+|v_k(x)|+1)\}.
\]
Since $u_k-v_k\to 0$ in $L^p$, we see that 
\be\label{lemuc-1}
	\int_{A_k(\delta)} (|u_k(x)|+|v_k(x)|+1)^p\,dx \to 0~\text{as $k\to\infty$, for every $\delta$}.
\ee
In addition, $h_\infty$ is uniformly continuous on the compact set $\overline{O}\times \overline{B(0,1)}\subset \R^n\times \R^m$, with a modulus of continuity $\mu$, whence
\be\label{lemuc-2}
\begin{aligned}
	&\int_{\O\setminus A_k(\delta)} |h_\infty(x,u_k)-h_\infty(x,v_k)|\,dx\\
	&=\int_{\O\setminus A_k(\delta)} \left|h_\infty\Big(x,\frac{u_k}{|u_k|+|v_k|+1}\Big)
	-h_\infty\Big(x,\frac{v_k}{|u_k|+|v_k|+1}\Big)\right| (|u_k(x)|+|v_k(x)|+1)^p\,dx\\
	&\leq \int_{\O\setminus A_k(\delta)} \mu(\delta)(|u_k(x)|+|v_k(x)|+1)^p\,dx\\
	&\leq \mu(\delta) C \underset{\delta\to 0}{\longrightarrow} 0\quad\text{uniformly in $k$},
\end{aligned}
\ee
where we also used that $\{u_k\}$ and $\{v_k\}$ are bounded in $L^p$. Combining \eqref{lemuc-1} and \eqref{lemuc-2}, $\|h_\infty(\cdot,u_k(\cdot))-h_\infty(\cdot,v_k(\cdot))\|_{L^1}$ can be made arbitrarily small, first choosing $\delta$ small enough and then $k$ large, depending on $\delta$. 
\hfill
$\Box$

\bigskip\bigskip
\begin{minipage}[t]{14cm}
{\bf Acknowledgment:}
We acknowledge the support through the project CZ01-DE03/2013-2014/DAAD-56269992 (PPP program). Moreover, MK and GP   were partly supported by the grants GA\v{C}R  P201/12/0671  and GAUK 267310, respectively. 
\end{minipage}

\end{sloppypar}

\end{document}